\documentclass[11pt]{article}
\usepackage{amssymb}
\usepackage{amsbsy}
\usepackage{cite}
\usepackage{hyperref}
\usepackage[numbers,sort&compress]{natbib}
\begin{document}
\title{\Large \bf The even $L_p$ Gaussian dual Minkowski problem
\thanks{Research is supported by the Natural Science
Foundation of China (Grant No.12161078)}}
\author{\small \bf Wei Shi$^{1}$, {Jiancheng Liu$^{2}$\thanks{Corresponding author, E-mail address: liujc@nwnu.edu.cn.}}
\\ \small  (College of Mathematics and Statistics, Northwest Normal University, Lanzhou, 730070, China)}
\date{}
\maketitle
\vskip 20pt

\begin{center}
\begin{minipage}{12cm}
\small
 {\bf Abstract:}
 The even Gaussian dual Minkowski problem studied by Feng, Hu and Xu. In this paper, we consider the even $L_p$ Gaussian dual Minkowski problem for $p>1$. The existence of $o$-symmetric solution in the case $p>1$ is obtained.



 \vskip 0.1cm
\end{minipage}
\end{center}

\vskip 20pt

\section{\bf Introduction}
\ \ \ \ The classical Brunn-Minkowski theory, also known as the theory of mixed volumes, is the core of convex geometric analysis. It originated with Minkowski when he combined his concept of mixed volume with the Brunn-Minkowski inequality. One of the Minkowski's major contributions towards the theory was to show how his theory could be developed from a few basic concepts, such as support function, vector sum, volumes and quermassintegrals. Gardner's and Schneider's classical treatises (\cite{gard,sch}) are two excellent surveys and source of references.

The center of the theory is the study of geometric invariant such as quermassintegrals, including volume, surface area, mean width and so on, and geometric measures such as area measures and Federer’s curvature measures. Geometric measures are closely related to geometric invariants and are differentials of the invariants when viewed as geometric functionals of convex bodies. The most significant functional of convex body in Euclidean space $\mathbb{R}^n$ is volume which is denoted by $V$. The unit sphere in $\mathbb{R}^n$ is denoted by $\mathbb{S}^{n-1}$. A convex body $K$ in $\mathbb{R}^n$ is a compact convex set with non-empty interior. Denote by $\mathcal{K}_o^n$ the set of all convex bodies that contain the origin in their interiors in $\mathbb{R}^n$ and by $\mathcal{K}_e^n$ the set of all origin-symmetric convex bodies in $\mathbb{R}^n$. The standard inner product of the vectors $x, y\in\mathbb{R}^n$ is denoted by $\langle x,y\rangle$. We write $|x|=\sqrt{\langle x,x\rangle}$, and use $\mathcal{H}^k$ to denote the $k$-dimensional Hausdorff measure and $\partial K$ the boundary of a convex body $K$.

Aleksandrov(\cite{alek})first established the following variational formula for volumes
$$\lim_{t\rightarrow 0}\frac{V(K+tL)-V(K)}{t}= \int_{S^{n-1}}h_L(v)dS_K(v),$$
where $K+tL=\{x+ty: x\in K \ \mathrm{and} \ y\in L\}$, $h_L:\mathbb{S}^{n-1}\rightarrow \mathbb{R}$ is the support function of convex body $L$. It naturally generates the surface area measure $S_K$ of convex body $K$. The surface area measure is the most studied measure in the classical Brunn-Minkowski theory, and the corresponding problem of characterizing the surface area measure is the classical Minkowski problem, which reads: {\it given a finite Borel measure $\mu$ on $\mathbb{S}^{n-1}$, what are the necessary and sufficient conditions on $\mu$ does there exist a convex body $K$ in $\mathbb{R}^n$ such that surface area measure $S_K(\cdot)=\mu(\cdot)$? If $K$ exists, is it unique?} The classical Minkowski problem was solved by Minkowski(\cite{MK})himself in the polytope case, and later its complete solution was given by Aleksandrov(\cite{alek}), Fenchel and Jessen(\cite{fen}). Solving the Minkowski problem is equivalent to solving the Monge-Ampère equation in the smooth case. See, e.g., Nirenberg(\cite{NL}), Pogorelov(\cite{PG}),Cheng and Yau(\cite{cy}) and Caffarelli(\cite{caff}).

The $L_p$-Brunn-Minkowski theory is a crucial extension of the classical Brunn-Minkowski theory whose roots date back to the Firey’s $L_p$-combination in 1962(\cite{fy}). But it was not actively developed until the emergence of the concept of the $L_p$-surface area measure in Lutwak (\cite{lu1}) in the 1990s. The $L_p$-Minkowski problem aims to characterize $L_p$-surface area measure defined as $S_{p,K}=h_K^{1-p}S_K$ for a convex body $K$. Clearly, when $p=1$, the $L_p$-surface area measure is just the classical surface area measure, and hence the $L_1$
Minkowski problem is the classical Minkowski problem. See \cite{alek1,MK1,TW}. Another notable case is $p=0$, where $S_0(K,\cdot)$ is also known as the cone volume measure, and the corresponding Minkowski problem is the logarithmic Minkowski problem. This case is largely unsolved, and results can be found in \cite{bk,bk1,cl,st,zg}. The $p=-n$ case (also largely open) is called the centro-affine Minkowski problem posed by Chou and Wang\cite{cw}. See [\cite{jian2016mirror,zg11}] for progress on the problem. For the $p>1$ case, Lutwak(\cite{lu1}) showed existence for even measures, and then Chou and Wang\cite{cw} proved existence for general measures. More results on this can be found in\cite{cl1,lw,cwx,hug}. Other important and related areas of research in $L_p$ Brunn-Minkowski theory include $L_p$-affine isoperimetric inequalities and sharp affine Sobolev inequalities. For results on this, see\cite{ci,hab,hab1,kn,lyz,lyz1,lyz2,zgy}.

In the 1970s, Lutwak\cite{luw} introduced the dual Brunn-Minkowski theory, establishing it as the dual counterpart to the classical Brunn-Minkowski theory. Within this framework, the classical quermassintegrals have corresponding dual quermassintegrals. While numerous dualities exist between the two theories, a notable gap has been the absence of a duality between the area measure and the Federer's curvature measure. This challenge was effectively addressed when Huang, Lutwak, Yang and Zhang\cite{huang} made significant advancements by discovering new geometric measures called dual curvature measures $\widetilde{C}_q(K,\cdot)$, which operate as duals to the Federer's curvature measure. In dual Brunn-Minkowski theory, significant research centers on the Minkowski problem. Naturally, Huang, Lutwak, Yang and Zhang\cite{huang} put forward the Minkowski problem concerning the dual curvature measures $\widetilde{C}_q(K,\cdot)$, which called the dual Minkowski problem, it asks: {\it for $q\in\mathbb{R}$ and a given non-zero finite Borel measure $\mu$ on $\mathbb{S}^{n-1}$, does there exist a convex body $K$ such that $\mu=\widetilde{C}_q(K,\cdot)$?} When $q=0$,the dual Minkowski problem becomes the classical Aleksandrov problem completely solved by Aleksandrov\cite{alek2} himself using a topological argument. For research of the dual Minkowski problem, please refer to the references\cite{bhh,blyzz,hep,huangz,zhao,zhao1}.

In recent years, the Minkowski problems have been gradually extended from the $L_p$ space to some other spaces, such as Orlicz space and the Gaussian probability space and so on, see e.g., \cite{ghwxy,ghxy,hlyz,hqz,jl,liqi,wuyu,zbc,laj,lyz11,lyz22,zx,ghw,zg1}. Huang, Xi, and Zhao\cite{hxz} studied the Minkowski problem in Gaussian probability space that is the natural analogue to the one for the Minkowski problem in $R^n$ which the volume and the surface area measure are replaced by Gaussian volume and the first variation of Gaussian volume respectively. The Gaussian Minkowski problem was posed in \cite{hxz}:{\it given a finite Borel measure $\mu$ on $\mathbb{S}^{n-1}$, what are the necessary and sufficient conditions on $\mu$ so that there exists a convex body $K$ whose Gaussian surface area measure satisfies $\mu=S_{\gamma_n, K}$? If $K$ exists, is it unique?} Regarding the Minkowski problem in Gaussian probability spaces, many scholars have already done excellent research in this regard. For specific details, please refer to references \cite{chlz,feng1,feng2,hu,ljq,sx,tang}.

The Gauss dual quermassintegral $\widetilde{V}_{\gamma_n,q}(K)$ of $K\in \mathcal{K}_o^n$ for $q>0$ is recently introduced by Feng, Hu and Xu\cite{feng3},which can be stated as follows
$$\widetilde{V}_{\gamma_n,q}(K)=\int_K e^{-\frac{|x|^2}{2}}|x|^{q-n}d\mathcal{H}^{n-1}(x), \eqno(1.1)$$
based on the definition of the Gaussian dual quermassintegrals mentioned above, the Gaussian dual quermassintegrals was obtained by using the variational formula in \cite{feng3}:
$$\lim_{t\rightarrow 0}\frac{\widetilde{V}_{\gamma_n,q}([h_t])-\widetilde{V}_{\gamma_n,q}(K)}{t}=\int_{S^{n-1}}f(v)d\widetilde{C}_{\gamma_n,q}(K,v),$$
where
$$h_t=e^{tf(v)}h_K(v)$$
for a continuous function $f: \mathbb{S}^{n-1}\rightarrow \mathbb{R}$ and
$$[h_t]=\{x\in \mathbb{R}^n:\langle x,v\rangle\leq h_t(v) \ \textmd{for all }v\in \mathbb{S}^{n-1}\}$$
is the Wulff shape of $h_t: \mathbb{S}^{n-1}\rightarrow(0,\infty)$.

For $q>0$ and each Borel measurable set $\eta\subset \mathbb{S}^{n-1}$, the Gaussian dual curvature measure is given by $$\widetilde{C}_{\gamma_n,q}(K,\eta)=\int_{\nu^{-1}_K(\eta)}\langle x,\nu_K(x)\rangle e^{-\frac{|x|^2}{2}}|x|^{q-n}d\mathcal{H}^{n-1}(x),$$
where $\nu^{-1}_K: \mathbb{S}^{n-1}\rightarrow\partial K$ is the inverse Gauss map.

In this article, based on the research in the above-mentioned \cite{feng3}, we constructed the $L_p$ Gaussian dual curvature measures $\widetilde{C}_{p,\gamma_n,q}(K,v)$ and obtained the following $L_p$ variational formula:
$$\lim_{t\rightarrow 0}\frac{\widetilde{V}_{\gamma_n,q}([h_t])-\widetilde{V}_{\gamma_n,q}(K)}{t}=\frac{1}{p}\int_{S^{n-1}}f(v)^p d\widetilde{C}_{p,\gamma_n,q}(K,v),$$
Note that here $q>0$ and $h_t(v)=(h_K(v)^p+tf(v)^p)^{\frac{1}{p}}$ for $p>0$.
From the above $L_p$ variational formula, we can naturally derive the definition of the $L_p$ dual Gaussian curvature measure: for any convex body $K$ and Borel measurable set $\eta\subset \mathbb{S}^{n-1}$, we have
$$\widetilde{C}_{p,\gamma_n,q}(K,v)=\int_{\nu_K^{-1}(\eta)}\langle x,\nu_K(x)\rangle^{1-p} |x|^{q-n} e^{-\frac{|x|^2}{2}}d\mathcal{H}^{n-1}(x).$$
{ \bf The $L_p$ Gaussian dual Minkowski problem} reads: {\it For $p,q>1$ and a given non-zero finite Borel measure $\mu$ on the unit sphere $\mathbb{S}^{n-1}$, does there exist a convex body $K\in \mathcal{K}_o^n$ such that
$$\mu=\widetilde{C}_{p,\gamma_n,q}(K,\cdot)?$$}
When the given measure $\mu$ has a density function $f$, then the $L_p$ Gaussian dual Minkowski problem is equivalent to solving the following Monge-Ampère equation on $\mathbb{S}^{n-1}$,
$$e^{-\frac{|\nabla h_K|^2+h_K^2}{2}}h_K^{-p}h_K(|\nabla h_K|^2+h_K^2)^{\frac{q-n}{2}}\textmd{det}(\nabla^2 h_K+h_K I)=f,$$
where $\nabla h_K$ and $\nabla^2 h_K$ represent the gradient and the Hessian matrix of the support function $h_K$ with respect to an orthogonal basis on $\mathbb{S}^{n-1}$ respectively, and $I$ is the identity matrix.

The main purpose of this paper is to solve the existence of solutions to the $L_p$ Gaussian dual Minkowski problem in the case of $p,q>1$. The Aleksandrov variational method in \cite{feng3} is mainly adopted. Now, the theorem is stated as follows.

\vskip 0.05cm \noindent{\bf Theorem 1.1}~~{\it For $p,q>1$, let $\mu$ be a nonzero finite Borel measure on $\mathbb{S}^{n-1}$, and there is an even density function $f$ on $\mathbb{S}^{n-1}$ with $d\mu=fdv$. Then there exists an origin-symmetric convex body $K\in\mathcal{K}_e^n$ such that
$$\frac{\mu}{|\mu|}=\frac{\widetilde{C}_{p,\gamma_n,q}(K,\cdot)}{|\widetilde{C}_{p,\gamma_n,q}(K,\cdot)|}.$$}

This paper is organized as follows. In section 2, we mainly present fundamental concepts related to convex bodies. Section 3 discusses properties of $L_p$ Gaussian dual curvature measures. In section 4, we focus on establishing the variational formula for $L_p$ Gaussian dual quermassintegrals, which is essential for deriving the $L_p$ Gaussian dual Minkowski problem when $p$ is greater than zero. The paper concludes with a proof of the main theorem.

\section{\bf Preliminaries}
\ \ \ \
In this section we introduce notation and collect basic facts from classical theory of convex bodies that we use in the paper. For more details, see \cite{sch}.

Let $\mathbb{S}^{n-1}$ be the $n$-dimensional unit sphere, thus, $S^{n-1}=\{x\in\mathbb{R}^n:|x|=1\}$. The set of all continuous functions defined on $\mathbb{S}^{n-1}$ is denoted by $C(\mathbb{S}^{n-1})$ and will always be viewed as equipped with the max-norm metric:
$$\|f-g\|_\infty =\max_{u\in \mathbb{S}^{n-1}}|f(u)-g(u)|.$$
The set of all even functions in $C(\mathbb{S}^{n-1})$ will write $C_e(\mathbb{S}^{n-1})$. Let $C^+(\mathbb{S}^{n-1})$ be the set of all positive functions in $C(\mathbb{S}^{n-1})$. $C_e^+(S^{n-1})$ for the set all even and positive functions in $C(\mathbb{S}^{n-1})$.

The support function $h_K$ of a compact convex subset $K$ in $\mathbb{R}^n$ is defined by
$$h_K(x)=\max\{\langle x,y\rangle: y\in K\}$$
for any $x\in \mathbb{R}^n$. Note that the support function $h_K$ is a continuous function and homogeneous of degree 1.

Define the radial function, $\rho_K=\rho(K,\cdot):\mathbb{R}^n\setminus\{0\}\rightarrow[0,+\infty)$, of a compact star-shaped (about the origin) set $K$ in $\mathbb{R}^n$ by (see \cite{sch})
$$\rho(K,x)=\max\{\lambda\geq0:\lambda x\in K\}, \ \ \ x\in\mathbb{R}^n\setminus\{0\}.$$
The radial function is a continuous function and homogeneous of degree -1. A direct observation is that $\partial K=\{\rho_K(u)u: u\in S^{n-1}\}$.

For each $f\in C^+(\mathbb{S}^{n-1})$, the{\it Wulff shape} $[f]$ generated by the function $f$ is a convex body defined by
$$[f]=\bigcap_{v\in \mathbb{S}^{n-1}}\{x\in \mathbb{R}^n: \langle x,v\rangle \leq h(v)\}.$$
It is easy to see that for $f\in C^+(\mathbb{S}^{n-1})$,
$$h_{[f]}\leq f, \eqno(2.1)$$
and for $K\in\mathcal{K}_o^n$,
$$[h_K]=K. \eqno(2.2)$$

For $K, L\subset\mathcal{K}^n$ and real $a, b>0$, the Minkowski combination $aK+bL\in \mathcal{K}^n$ is defined as follows
$$aK+bL=\{ax+by: x\in K \ \textmd{and} \ y\in L\},$$
and its support function is
$$h_{aK+bL}=ah_K+bh_L.$$

The {\it $L_p$ Minkowski combination} is the basic concept in the $L_p$-Brunn-Minkowski theory. Fix a real $p$. For $K, L\in\mathcal{K}_o^n$, and $a, b\geq 0$, define the $L_p$ Minkowski combination, $a\cdot K+_{p} b\cdot L\in \mathcal{K}_o^n$, via the Wulff shape:
$$a\cdot K+_{p} b\cdot L=[(ah_K^p(v)+bh_L^p(v))^\frac{1}{p}], \eqno(2.3)$$
when $p\neq 0$. Note that the notion of Wulff shape allows us to consider a $L_p$ Minkowski combination where either $a$ or $b$ may be negative, as long the function $ah_K^p+bh_L^p$ is strictly positive on $\mathbb{S}^{n-1}$. When $p=0$, define $a\cdot K+_0 b\cdot L$ via the Wulff shape
$a\cdot K+_0 b\cdot L=[h_K^{a}h_L^{b}]$.

We say that a sequence of convex bodies $\{K_i\}$ converges to a compact convex set $K$ in $\mathbb{R}^n$ if
$$\sup\{|h_{K_i}(v)-h_K(v)|: v\in\mathbb{S}^{n-1}\}\rightarrow 0,$$
as $i\rightarrow \infty.$

Define the supporting hyperplane, with outer unit normal $v\in \mathbb{S}^{n-1}$, tangent to $K\in\mathcal{K}_o^n$ by
$$H_K(v)=\{x\in \mathbb{R}^n: \langle x,v\rangle =h_K(v)\}.$$

For $\sigma\subset\partial K$, the spherical image, $\boldsymbol{\nu}_K(\sigma)$, of $\sigma$ is defined by
$$\boldsymbol{\nu}_K(\sigma)=\{v\in\mathbb{S}^{n-1}: x\in H_K(v) \ \mathrm{for} \ \mathrm{some} \ x\in \sigma\}\subset\mathbb{S}^{n-1}.$$

Suppose that $\sigma_K\subset\partial K$ is the set consisting of all $x\in\partial K$ for which the set $\boldsymbol{\nu}_K(\{x\})$, often abbreviated by $\boldsymbol{\nu}_K(x)$, contains more than a single element. It is well known that $\mathcal{H}^{n-1}(\sigma_K)=0$ (see Schneider \cite{sch}). The spherical image map is defined as $\nu_K: \partial K\backslash \sigma_K\rightarrow\mathbb{S}^{n-1}$. It will occasionally be convenient to abbreviate $\partial K\backslash \sigma_K$ by $\partial^{\prime} K$.

For $K\in\mathcal{K}_o^n$, define the radial map of $K$,
$$r_K: \mathbb{S}^{n-1}\rightarrow\partial K \ \mathrm{by} \ r_K(u)=\rho_K(u)u\in\partial K,$$
for $u\in\mathbb{S}^{n-1}$.

For $\omega\subset\mathbb{S}^{n-1}$, define the radial Gauss image of $\omega$ by
$$\boldsymbol{\alpha}_K(\omega)=\boldsymbol{\nu}_K (r_K(\omega))\subset \mathbb{S}^{n-1}.$$
Thus, for $u\in\mathbb{S}^{n-1}$,
$$\boldsymbol{\alpha}_K(u)=\{v\in\mathbb{S}^{n-1}: r_{K}(u)\in H_K(v)\}.$$

Define the radial Gauss map of the convex body $K\in \mathcal{K}_o^n$
$$\alpha_K: \mathbb{S}^{n-1}\backslash \omega_{K}\rightarrow \mathbb{S}^{n-1} \ \mathrm{by} \ \alpha_K=\nu_{K}\circ r_K,$$
where $\omega_K=\overline{\sigma_{K}}=r_K^{-1}(\sigma_K)$. Since $r_K^{-1}$ is a bi-Lipschitz map between the spaces $\partial K$ and $\mathbb{S}^{n-1}$ it follows that $\omega_K$ has spherical Lebesgue measure 0. Observe that if $u\in\mathbb{S}^{n-1}\backslash \omega_{K}$, then $\boldsymbol{\alpha}_K(u)$ contains only the element ${\alpha}_K(u)$. Note that since both $\nu_K$ and $r_K$ are continuous, $\alpha_K$ is continuous.

For any $\mathcal{H}^{n-1}-$integrable function $g: \partial K\rightarrow \mathbb{R}$, the following identity holds
$$\int_{\partial K}g(x)d\mathcal{H}^{n-1}(x)=\int_{S^{n-1}}g(\rho_K(u)u)F(u)du,$$
where $F$ is defined as $\mathcal{H}^{n-1}-$a.e. on $\mathbb{S}^{n-1}$ by
$$F(u)=\frac{(\rho_K(u))^n}{h_K(\alpha_K(u))}.$$

\section{\bf $L_p$ Gaussian dual curvature measure}
\ \ \ \ \ \
To give the definition of the $L_p$ Gaussian dual curvature measure, which is essentially the derivative of the Gaussian dual quermassintegrals with respect to the $L_p$ Minkowski combination, for this purpose, the following variational lemma on the support function of the $L_p$ Minkowski combination is what we need.

\vskip 0.05cm \noindent{\bf Lemma 3.1(\cite{ljq})}~~{\it For $p\neq 0$, let $K\in \mathcal{K}_o^n$, and $f: \mathbb{S}^{n-1}\rightarrow \mathbb{R}$ be a continuous function. For small enough $\delta>0$, and each $t\in(-\delta,\delta)$, we define the continuous function $h_t: \mathbb{S}^{n-1}\rightarrow(0,\infty)$ as
$$h_t(v)=(h_K(v)^p+tf(v)^p)^{\frac{1}{p}}, \  v\in S^{n-1}. \eqno(3.1)$$
Then,
$$\lim_{t\to0}\frac{\rho_{[h_t]}(u)-\rho_K(u)}{t}=\frac{f(\alpha_K(u))^p}{ph_K(\alpha_K(u))^p}\rho_K(u) \eqno(3.2)$$
holds for almost all $u\in \mathbb{S}^{n-1}$. In addition, there exists $M>0$ such that
$$|\rho_{[h_{t}]}(u)-\rho_{K}(u)|\leq M|t|, \eqno(3.3)$$
for all $u\in \mathbb{S}^{n-1}$ and $t\in(-\delta,\delta)$.}

Regarding the convergence of the Gaussian dual quermassintegrals, it has been demonstrated in detail in \cite{feng3}, and it plays a crucial role in proving the $L_p$ Gaussian dual Minkowski problem.

\vskip 0.05cm \noindent{\bf Lemma 3.2(\cite{feng3})}~~{\it For $q>1$, $K_i\in\mathcal{K}_o^n$ where $i=1,2,\cdots.$ If $K_i\rightarrow K_0\in\mathcal{K}_o^n$ with respect to the Hausdorff metric, then
$$\lim_{i\to\infty}\widetilde{V}_{\gamma_n,q}(K_i)=\widetilde{V}_{\gamma_n,q}(K_0). \eqno(3.4)$$}

Next, we show that the $L_p$ variational formula for the Gaussian dual quermassintegrals.

\vskip 0.05cm \noindent{\bf Theorem 3.1}~~{\it For $p,q>1$, let $K\in\mathcal{K}^n_o$ and $f: \mathbb{S}^{n-1}\rightarrow\mathbb{R}$ be a continuous function. For sufficiently small $\delta>0$ and each $t\in(-\delta,\delta)$, define $h_t=(h_K^p+tf^p)^{\frac{1}{p}}$. Then,
$$\lim_{t\rightarrow 0}\frac{\widetilde{V}_{\gamma_n,q}([h_t])-\widetilde{V}_{\gamma_n,q}(K)}{t}=\frac{1}{p}\int_{S^{n-1}}f(v)^p d\widetilde{C}_{p,\gamma_n,q}(K,v). \eqno(3.5)$$}
{\bf Proof.} Applying the polar coordinate formula, we have
$$\widetilde{V}_{\gamma_n,q}([h_t])=\int_{\mathbb{S}^{n-1}}\int_0^{\rho_{[h_t]}(u)}e^{-\frac{r^2}{2}}r^{q-1}drdu.$$
Noting that $h_t(v)=(h_K(v)^p+tf(v)^p)^{\frac{1}{p}}, \  v\in S^{n-1},$ and combining it with the (3.1), Taylor expansion of $\log(1+x)$, for sufficiently small $t>0$,
$$\log h_t=\log h_K+\frac{tf^p}{ph_K^p}+o(t,\cdot), \eqno(3.6)$$
where, the function $o(t,\cdot): \mathbb{S}^{n-1}\rightarrow \mathbb{R}$ is continuous and $\lim_{t\rightarrow 0}\frac{o(t,\cdot)}{t}=0$ uniformly on $\mathbb{S}^{n-1}$. Just as stated in \cite{huang}, we still call $[h_t]$ {\it a logarithmic family of Wulff shapes formed by $h_K$ and $\frac{f^p}{ph_K^p}$.}

Since $K\in \mathcal{K}_o^n$ and $h_t\rightarrow h_K$ uniformly as $t\rightarrow 0$, we get $[h_t]\rightarrow K$ by (2.2) and Aleksandrov’s convergence Lemma (\cite{sch}). Then radial function of $[h_t]$ converges to the radial function of $K$ correspondingly, i.e., $\rho_{[h_t]}\rightarrow \rho_K$, and there exist $c_0, c_1>0$ such that $\rho_{[h_t]}, \rho_K\in [c_0, c_1]$ for sufficiently small $t\in(-\delta,\delta)$.

For convenience, we let $G(s)=\int_{0}^s e^{-\frac{r^2}{2}}r^{q-1} dr$. By mean value theorem, there exists $\xi\in[c_0, c_1]$ such that
$$|G(\rho_{[h_t]}(u))-G(\rho_K(u))|=|G^\prime(\xi)||\rho_{[h_t]}(u)-\rho_K(u)|\leq M|G^\prime(\xi)||t|,$$
where $M$ is just as the positive number in (3.3). Since a continuous function on a closed interval $[c_0, c_1]$ must have maximum and minimum values, then, we have,
$$|G^\prime(\xi)|=|e^{-\frac{\xi^2}{2}}\xi^{q-1}|\leq M_1$$
for some constant $M_1>0$ when $\xi\in [c_0,c_1]$ by the definition of $G$. Therefore,
$$|G(\rho_{[h_t]}(u))-G(\rho_K(u))|\leq MM_1 |t|.$$

From formula (3.1), (3.2) and the definition (1.1) of the Gaussian dual quermassintegrals, employing the dominated convergence theorem, it follows that
$$\lim_{t\rightarrow 0}\frac{\widetilde{V}_{\gamma_n,q}([h_t])-\widetilde{V}_{\gamma_n,q}(K)}{t}=\lim_{t\rightarrow 0}\frac{\widetilde{V}_{\gamma_n,q}([(h_K(v)^p+tf(v)^p)^{\frac{1}{p}}])-\widetilde{V}_{\gamma_n,q}(K)}{t} \ \ \ \ \ \ \ \ \ \ \ \ \ \ \ \ \ \ \ \ \ \ \ $$
$$=\lim_{t\rightarrow 0}\frac{1}{t}\int_{\mathbb{S}^{n-1}}\int_{\rho_K(u)}^{\rho_{[h_t]}(u)}e^{-\frac{r^2}{2}}r^{q-1}drdu.$$
$$ \ \ \ \ \ \ \ \ \ \ \ \ \ \ \ \ \ \ \ \ \ \ \ \ \ =\frac{1}{p}\int_{\partial^\prime K}\langle x,\nu_K(x)\rangle^{1-p} |x|^{q-n} e^{-\frac{|x|^2}{2}}f(\nu_K(x))^p d\mathcal{H}^{n-1}(x) \ \ \ $$
$$=\frac{1}{p}\int_{\mathbb{S}^{n-1}}f(v)^p d\widetilde{C}_{p,\gamma_n,q}(K,v).\ \ \ \ \ \ \ \ \ \ $$
\hfill ${\square}$

That is, through the $L_p$ variational formula (3.5) of Gauss dual quermassintegrals, we have derived the $L_p$ dual Gaussian curvature measure.

\vskip 0.05cm \noindent{\bf Definition 3.1}~~ For $p,q>1$ and $K\in\mathcal{K}_o^n$, $\eta\subset\mathbb{S}^{n-1}$ is a Borel set on $\mathbb{S}^{n-1}$. $L_p$ dual Gaussian curvature measure is defined as below,
$$\widetilde{C}_{p,\gamma_n,q}(K,\eta)=\int_{\nu_K^{-1}(\eta)}\langle x,\nu_K(x) \rangle^{1-p} |x|^{q-n} e^{-\frac{|x|^2}{2}}d\mathcal{H}^{n-1}(x),$$
which is a Borel measure on $\mathbb{S}^{n-1}$.

The Gaussian dual curvature measure $\widetilde{C}_{\gamma_n,q}(K,\cdot)$ (see \cite{feng3}) is weakly convergent and absolutely continuous with respect to the surface area measure $S(K,\cdot)$. Therefore, it is essential to demonstrate that $\widetilde{C}_{p,\gamma_n,q}(K,\cdot)$ is weakly convergent in relation to the Hausdorff metric and absolutely continuous with respect to the surface area measure $S(K,\cdot)$.

\section{\bf Even $L_p$ dual-Gaussian Minkowski problem and its solution}
\ \ \ \ \
In this part, we mainly adopt the variational method to prove the existence of solutions to the $L_p$ dual Gaussian Minkowski problem when $p,q>1$. To this end, we first consider a related optimization problem and prove that the solutions satisfying this optimization problem are the solutions that meet the $L_p$ dual Gaussian Minkowski problem. Secondly, we illustrate that solutions that satisfy this optimization problem exist.

\vskip 0.05cm \noindent{\bf 4.1 An associated optimization problem}

For any nonzero finite Borel measure $\mu$ on the sphere $\mathbb{S}^{n-1}$, the total mass is denoted by $|\mu|=\int_{S^{n-1}}d\mu(v).$ For $p>1$,  $g\in C_e^+(\mathbb{S}^{n-1})$, define the functional
$$\widetilde{\Phi}: C_e^{+}(\mathbb{S}^{n-1})\rightarrow \mathbb{R}$$ by
$$\widetilde{\Phi}(g)=-\frac{\int_{\mathbb{S}^{n-1}}g(v)^p d\mu(v)}{p|\mu|}. \eqno(4.1)$$

We are looking for a function $g_0$ where $\widetilde{\Phi}$ reaches its maximum point, that is, considering the following optimization problem
$$ \sup\bigg\{\widetilde{\Phi}(g): \widetilde{V}_{\gamma_n,q}([g])=|\mu| \ \mathrm{and} \ g\in C_e^+(\mathbb{S}^{n-1})\bigg\}=\widetilde{\Phi}(g_0). \eqno(4.2)$$
The search can be restricted to support functions of convex bodies in $\mathcal{K}_e^n$. To see this, note that for each $g\in C^{+}(\mathbb{S}^{n-1})$ satisfies (2.1), thus,
$$h_{[g]}\leq g,$$
and
$$[h_{[g]}]=[g].$$
Since $\mu\geq 0$, together with the above facts, it follows that
$$\widetilde{\Phi}(g)\leq \widetilde{\Phi}(h_{[g]}), \widetilde{V}_{\gamma_n,q}([g])=\widetilde{V}_{\gamma_n,q}([h_{[g]}]).$$

Thus, the above optimization problem (4.2) can be transformed into the following optimization problem
$$\sup\bigg\{\Phi(K): \widetilde{V}_{\gamma_n,q}(K)=|\mu| \ \mathrm{and} \ K\in\mathcal{K}_e^n\bigg\}=\Phi(K_0), \eqno(4.3)$$
here, for $p>1$, $K\in\mathcal{K}_o^n$, define the functional $\Phi: \mathcal{K}_e^n\rightarrow\mathbb{R}$ by
$$\Phi(K)=-\frac{\int_{\mathbb{S}^{n-1}}h_K(v)^pd\mu(v)}{p|\mu|}.$$ And we easily see that $\Phi(K)=\widetilde{\Phi}(h_K)$.

\vskip 0.05cm \noindent{\bf Theorem 4.1} {\it \ Let $p, q>1$ and $\mu$ be a non-zero finite Borel measure on $\mathbb{S}^{n-1}$. If $K_0\in \mathcal{K}_e^n$, $\widetilde{V}_{\gamma_n,q}(K_0)=|\mu|$ and it is the following optimization problem's solution
$$\sup\bigg\{\Phi(K): \widetilde{V}_{\gamma_n,q}(K)=|\mu| \ \mathrm{and} \ K\in\mathcal{K}_e^n\bigg\}, $$
then
$$\frac{\widetilde{C}_{p,\gamma_n,q}(K,\cdot)}{|\widetilde{C}_{p,\gamma_n,q}(K,\cdot)|}=\frac{\mu}{|\mu|}.$$}

{\bf Proof.} Let $K_0\in\mathcal{K}_e^n$ be a maximum for optimization problem (4.3), thus, $h_{K_0}$ is a solution of optimization problem (4.2), we have $\widetilde{V}_{\gamma_n,q}([h_{K_0}])=|\mu|, h_{K_0}\in C_e^+(\mathbb{S}^{n-1}),$ and
$$\widetilde{\Phi}(h_{K_0})=\sup\bigg\{\widetilde{\Phi}(g): \widetilde{V}_{\gamma_n,q}([g])=|\mu| \ \mathrm{and} \ g\in C_e^+(\mathbb{S}^{n-1})\bigg\}.$$

For any $g\in C_e^+(\mathbb{S}^{n-1})$, $\lambda\in\mathbb{R}$ and $t\in(-\delta,\delta)$ where $\delta>0$ is sufficiently small, define $h_t\in C_e^+(\mathbb{S}^{n-1})$ by
$$h_t=(h_{K_0}^p+tg^p)^{\frac{1}{p}},$$
the above formula is actually equivalent to (3.6). Now it is known that $h_{K_0}$ is the maximum point, so
$$\widetilde{\Phi}(h_{K_0})\geq \widetilde{\Phi}(h_{t}). \eqno(4.4)$$

Let $\Gamma(t,\lambda)=\widetilde{\Phi}(h_{t})+\lambda(\widetilde{V}_{\gamma_n,q}([h_t])-|\mu|)$,
putting $F(v,t)=h_{K_0}(v)^p+tg(v)^p$,    $G(v,t)=\int_0^{\rho_{[h_t]}(v)}e^{-\frac{r^2}{2}}r^{q-1}dr$, and $(v,t)\in \mathbb{S}^{n-1}\times(-\delta,\delta)$. Observe that both $F(v,t)$ and $G(v,t)$ are integrable as functions of $v$ on $\mathbb{S}^{n-1}$, and are differentiable as functions of $t$ on $(-\delta,\delta)$. Obviously, $|\frac{\partial F(v,t)}{\partial t}|=|g(v)^p|\leq L$ since $g(v)\in C(\mathbb{S}^{n-1})$, and
$$|\frac{\partial G(v,t)}{\partial t}|=e^{-\frac{{\rho_{[h_t]}(v)}^2}{2}}{\rho_{[h_t]}(v)}^{q-1}|\frac{\partial \rho_{[h_t]}(v)}{\partial t}|\leq M e^{-\frac{{\rho_{[h_t]}(v)}^2}{2}}{\rho_{[h_t]}(v)}^{q-1}$$
via (3.3). Here, this constant $L$ and $
M e^{-\frac{{\rho_{[h_t]}(v)}^2}{2}}{\rho_{[h_t]}(v)}^{q-1}$ are integrable as functions of $v$ on $\mathbb{S}^{n-1}$. So, utilizing differentiation under the integral sign, by (4.4), (3.5), (4.1), (3.6) and (3.2), we have
$$0=\frac{\partial \Gamma(t,\lambda)}{\partial t}\bigg|_{t=0} \ \ \ \ \ \ \ \ \ \ \ \ \ \ \ \ \ \ \ \ \ \ \ \ \ \ \ \ \ \ \ \ \ \ \ \ \ \ \ \ \ \ \ \ \ \ \ \ \ \ \ \ \ \ \ \ \ \ \ \ \ \ \ \ \ \ $$
$$=\frac{\partial}{\partial t}\bigg(\widetilde{\Phi}(h_{t})+\lambda(\widetilde{V}_{\gamma_n,q}([h_t])-|\mu|)\bigg)\bigg|_{t=0} \ \ \ \ \ \ \ \ \ \ \ \ \ \ \ \ \ \ \ \ \ \ \ \ \ \ \ \ \ \ \ \ $$
$$=\frac{\partial}{\partial t}\bigg(-\frac{\int_{\mathbb{S}^{n-1}}(h_{K_0}(v)^p+tg(v)^p)d\mu(v)}{p|\mu|}\bigg)\bigg|_{t=0} \ \ \ \ \ \ \ \ \ \ \ \ \ \ \ \ \ \ \ \ \ \ \ \ \ $$
$$+\lambda\frac{\partial}{\partial t}\bigg(\int_{\mathbb{S}^{n-1}}\int_0^{\rho_{[h_t]}(v)}e^{-\frac{r^2}{2}}r^{q-1}drdv-|\mu|\bigg)\bigg|_{t=0}  \ \ \ \ \ \ \ \ \ \ \ \ \ \ \ \ \ \ \ $$
$$ \ \ \ \ \ \ \ \ \ \ \ \ \ \ \ =\bigg(-\frac{1}{p|\mu|}\int_{S^{n-1}}g(v)^p d\mu(v)+\lambda\int_{S^{n-1}}e^{-\frac{{\rho_{[h_t]}(v)}^2}{2}}{\rho_{[h_t]}(v)}^{q-1}\frac{\partial \rho_{[h_t]}(v)}{\partial t}dv\bigg)\bigg|_{t=0}$$
$$ \ \ =\bigg(-\frac{1}{p|\mu|}\int_{S^{n-1}}g(v)^p d\mu(v)+\frac{\lambda}{p}\int_{\mathbb{S}^{n-1}}g(v)^p d\widetilde{C}_{p,\gamma_n,q}(K_0,v)\bigg)\bigg|_{t=0}, $$
that is
$$\frac{\int_{\mathbb{S}^{n-1}}g(v)^p d\mu(v)}{|\mu|}=\lambda\int_{\mathbb{S}^{n-1}}g(v)^p d\widetilde{C}_{p,\gamma_{n},q}(K_0,v). \eqno(4.5)$$
The above formula holds for any $g\in C^+(\mathbb{S}^{n-1})$, then
$$\mu=\lambda|\mu|\widetilde{C}_{p,\gamma_{n},q}(K_0,\cdot). \eqno(4.6)$$
Putting $g(v)=1$ in (4.5), then
$$\lambda=\frac{1}{|\widetilde{C}_{p,\gamma_{n},q}(K_0,\cdot)|},$$
together with (4.6), finally, it follows that
$$\frac{\mu}{|\mu|}=\frac{\widetilde{C}_{p,\gamma_{n},q}(K_0,\cdot)}{|\widetilde{C}_{p,\gamma_{n},q}(K_0,\cdot)|}.$$ \hfill${\square}$

\vskip 0.05cm \noindent{\bf 4.2 Existence of an optimizer}

\vskip 0.05cm \noindent{\bf Theorem 4.2} {\it \ Let $p,q>1$, and the measure $\mu$ is absolutely continuous with respect to the Lebesgue measure on $\mathbb{S}^{n-1}$ and has an even density $g(u)$, i.e., $d\mu(u)=g(u)du$, such that $\frac{1}{M_{0}}\leq g(u)\leq M_{0}$ for some constant $M_0>0$. Then there exists a convex body $K_{0}\in\mathcal{K}_e^n$ satisfies $\widetilde{V}_{\gamma_n,q}(K_0)=|\mu|$ and
$$\Phi(K_0)=\sup\bigg\{\Phi(K): \widetilde{V}_{\gamma_n,q}(K)=|\mu| \ \mathrm{and} \ K\in\mathcal{K}_e^n\bigg\}.$$}

{\bf Proof.} Suppose that $K_i\subset\mathcal{K}_o^n$ is a maximizing sequence to optimization problem (4.3), thus, it follows that
$$\widetilde{V}_{\gamma_n,q}(K_i)=|\mu|,$$
and
$$\lim_{i\rightarrow\infty}\Phi(K_i)=\sup\bigg\{\Phi(K): \widetilde{V}_{\gamma_n,q}(K)=|\mu| \ \mathrm{and} \ K\in\mathcal{K}_e^n\bigg\}. \eqno(4.7)$$

We claim that $K_i$ is uniformly bounded. Let $2R_i$ denote the length of the longest segment in $K_i$. Since the $K_i$ are origin-symmetric, there exist $v_i\in\mathbb{S}^{n-1}$ such that $R_i v_i$, $-R_i v_i$ are the endpoints of this segment. Since this segment is contained in $K_i$, its support function is dominated by the support function of $K_i$, thus, $R_i|\langle v_i, u\rangle|\leq h_{K_i}(u)$, for all $u\in\mathbb{S}^{n-1}$. Based on this fact and combined with $\frac{1}{M_{0}}\leq g(u)\leq M_{0}$ for some constant $M_0>0$, we have
$$\Phi(K_i)=\widetilde{\Phi}(h_{K_i})=-\frac{\int_{\mathbb{S}^{n-1}}h_{K_i}(u)^p d\mu(u)}{p|\mu|}$$
$$ \ \ \ \ \ \ \ \ \ \ \ \ \ \ \ \ \ \ \ \ \ \ \ \ \ \ \ \ \leq -\frac{{R_i}^p\int_{\mathbb{S}^{n-1}}|\langle v_i, u\rangle|^p d\mu(u)}{p|\mu|}$$
$$ \ \ \ \ \ \ \ \ \ \ \ \ \ \ \ \ \ \ \ \ \ \ \ \ \ \ \ \ \ \ \ \ =-\frac{{R_i}^p\int_{\mathbb{S}^{n-1}}|\langle v_i, u\rangle|^p g(u)d(u)}{p|\mu|}$$
$$ \ \ \ \ \ \ \ \ \ \ \ \ \ \ \ \ \ \ \ \ \ \ \ \ \ \ \ =-\frac{{R_i}^p\int_{\mathbb{S}^{n-1}}|\langle v_i, u\rangle|^p d(u)}{pM_0|\mu|}.$$

Obviously, $h_{K_i}$ has a uniformly positive upper bound. If not, $\lim_{i\rightarrow\infty}R_i=\infty,$ together with above formula and note that $|\langle v_i, u\rangle|^p$  is integrable on $\mathbb{S}^{n-1}$, it follows that
$$\Phi(K_i)\leq -\frac{{R_i}^p\int_{\mathbb{S}^{n-1}}|\langle v_i, u\rangle|^p d(u)}{pM_0|\mu|}\rightarrow -\infty,$$
but this contradicts (4.7). So we conclude that $h_{K_i}$ has a uniformly positive upper bound.

Next, we will show that $h_{K_i}$ has a uniformly positive lower bound. In fact, this has been illustrated in \cite{feng3}, that is, there exists a positive constant $M_1$,  such that
$$\min_{v\in S^{n-1}}h_{K_i}(v)\geq M_1$$
holds.

Combining the above discussions, it can be finally concluded that there exists a positive constant $M_2$ independent of i such that
$$\frac{1}{M_2}\leq h_{K_i}\leq M_2.$$
By Blaschke selection theorem in \cite{sch}, there is a convergent subsequence of $K_i$, still denoted by $K_i$, converges to an origin-symmetric convex body $K_0$ of $\mathbb{R}^{n}$. By the continuity of $\Phi$ and (3.4), then
$$|\mu|=\lim_{i\rightarrow \infty}\widetilde{V}_{\gamma_n,q}(K_i)=\widetilde{V}_{\gamma_n,q}(K_0),$$
and
$$\lim_{i\rightarrow\infty}\Phi(K_i)=\sup\bigg\{\Phi(K): \widetilde{V}_{\gamma_n,q}(K)=|\mu| \ \mathrm{and} \ K\in\mathcal{K}_e^n\bigg\}=\Phi(K_0).$$

\vskip 0.05cm \noindent{\bf 4.3 Existence of solution to the even $L_p$ Gaussian dual Minkowski problem}

Combining Theorem 4.1 and Theorem 4.2, we get the main existence Theorem 1.1 for the even $L_p$ Gaussian dual Minkowski problem stated in the introduction.

\vskip 0.05cm \noindent{\bf Theorem 4.3}~~{\it For $p,q>0$, let $\mu$ be a nonzero finite Borel measure on $\mathbb{S}^{n-1}$, and there is an even density function $f$ on $\mathbb{S}^{n-1}$ with $d\mu=fdv$. Then there exists an origin-symmetric convex body $K\in\mathcal{K}_e^n$ such that
$$\frac{\mu}{|\mu|}=\frac{\widetilde{C}_{p,\gamma_n,q}(K,\cdot)}{|\widetilde{C}_{p,\gamma_n,q}(K,\cdot)|}.$$}

  \vskip 1.0cm
\bibliographystyle{abbrv}
\bibliography{name}

\clearpage

\end{document}